# TOTAL DILATIONS

Jean-Christophe Bourin

Les Coteaux, 8 rue Henri Durel 78510 Triel, France

E-mail: bourinjc@club-internet.fr

**abstract**

(1) Let $A$ be an operator on a space $\mathcal{H}$ of even finite dimension. Then for some decomposition $\mathcal{H} = \mathcal{F} \oplus \mathcal{F}^\perp$, the compressions of $A$ onto $\mathcal{F}$ and $\mathcal{F}^\perp$ are unitarily equivalent. (2) Let $\{A_j\}_{j=0}^n$ be a family of strictly positive operators on a space $\mathcal{H}$. Then, for some integer $k$, we can dilate each $A_j$ into a positive operator $B_j$ on $\oplus^k \mathcal{H}$ in such a way that: (i) The operator diagonal of $B_j$ consists of a repetition of $A_j$. (ii) There exist a positive operator $B$ on $\oplus^k \mathcal{H}$ and an increasing function $f_j : (0, \infty) \longrightarrow (0, \infty)$ such that $B_j = f_j(B)$.

**keywords** dilation, positive operators

**AMS subjects classification** 47A20

**0. Introduction**

This paper is a continuation of a subsection of [2] entitled "commuting dilations". We recall our definitions and notations. A pair of positive (semi-definite) operators $A$ and $B$ on a finite dimensional Hilbert space $\mathcal{H}$, $\dim \mathcal{H} = d$, is said to be a monotone pair of positive operators, or a positive monotone pair, if there exists an orthonormal basis $\{e_k\}_{k=1}^d$ such that

$$A = \sum_{k=1}^d \mu_k(A) e_k \otimes e_k \qquad \text{and} \qquad B = \sum_{k=1}^d \mu_k(B) e_k \otimes e_k$$

where the numbers $\mu_k(\cdot)$ are the singular values arranged in decreasing order and counted with their multiplicities and $e_k \otimes e_k$ is the rank one projection associated with $e_k$. On the other hand, if

$$A = \sum_{k=1}^d \mu_k(A) e_k \otimes e_k \qquad \text{and} \qquad B = \sum_{k=1}^d \mu_{d+1-k}(B) e_k \otimes e_k$$

we say that $(A, B)$ is an antimonotone pair of positive operators. It is easy to define the notion of a monotone family $\{A_j\}_{j=0}^n$ of positive operators. Furthermore, this notion can be extended to the notion of a monotone family of hermitian operators $\{A_j\}_{j=0}^n$ by requiring that there is a (hilbertian) basis $\{e_k\}$ for which

$$A_j = \sum_{k \geq 1} \lambda_k(A_j) e_k \otimes e_k, \qquad 0 \leq j \leq n,$$

where $\lambda_k(\cdot)$ are the eigenvalues arranged in decreasing order and counted with their multiplicities. Setting $A = \sum_{k=1}^d (d-k) e_k \otimes e_k$, we note that $A_j = f_j(A)$ for some increasing functions $f_j$.

Positive, monotone pairs $(A, B)$ well behave in respect to the compression to a subspace $\mathcal{E}$ of $\mathcal{H}$ (we recall this classical notion in Section 1). For instance we proved [2, Corollaries 2.2 and 2.3] that

$$\lambda_k(A_\mathcal{E} B_\mathcal{E}) \leq \lambda_k((AB)_\mathcal{E})$$

and

$$\lambda_k(A_\mathcal{E} B_\mathcal{E} A_\mathcal{E}) \leq \lambda_k((ABA)_\mathcal{E})$$

for all $k$. From the first inequality we derived

$$\det A_\mathcal{E} \cdot \det B_\mathcal{E} \leq \det(AB)_\mathcal{E}$$

while we showed that, in case of an antimonotone pair $(A, B)$ and a hyperplane $\mathcal{E}$, we have the opposite inequality

$$\det A_\mathcal{E} \cdot \det B_\mathcal{E} \geq \det(AB)_\mathcal{E}.$$

These results suggest the following question: Given a finite family of positive operators, how can we dilate them into a positive, monotone family? This paper precisely deals with the construction of such monotone dilations. However it appears that the dilations built up have the additional property to be *total* dilations. This notion is discussed in Section 1; the main result herein is the proof of the following fact:

*Any $2n$-by-$2n$ matrix $A$ is unitarily equivalent to a matrix of the form*

$$\begin{pmatrix} B & \star \\ \star & B \end{pmatrix}$$

*in which $B$ is some $n$-by-$n$ matrix and the stars hold for unspecified entries.*

We devote Section 2 to the study of monotone dilations. This section is divided in two subsections; the first one presents results whose proofs have an algorithmic nature while the second one gives more theoretical facts.

**1. Dilations and total dilations.**

Let $B$ be an operator on a space $\mathcal{H}$ and let $\mathcal{E}$ be a subspace of $\mathcal{H}$. Denote by $E$ the projection onto $\mathcal{E}$. The restriction of $EB$ to $\mathcal{E}$, denoted by $B_{\mathcal{E}}$, is the compression of $B$ to $\mathcal{E}$. Therefore, in respect to the decomposition $\mathcal{H} = \mathcal{E} \oplus \mathcal{E}^{\perp}$, we may write

$$B = \begin{pmatrix} B_{\mathcal{E}} & \star \\ \star & \star \end{pmatrix}.$$

The notion of compression has a natural extension: If $A$ is an operator on a space $\mathcal{F}$ with $\dim \mathcal{F} \leq \dim \mathcal{H}$, we still say that $A$ is a compression of $B$ if there is an isometry $V : \mathcal{F} \longrightarrow \mathcal{H}$ such that $A = V^*BV$. Thus, identifying $A$ with $VAV^*$ (equivalently, identifiyng $\mathcal{F}$ and $V(\mathcal{F})$), we can write

$$B = \begin{pmatrix} A & \star \\ \star & \star \end{pmatrix}.$$

One also says that $B$ dilates $A$ or that $B$ is a dilation of $A$.

Denote by $\oplus^k \mathcal{H}$ the direct sum $\mathcal{H} \oplus \ldots \oplus \mathcal{H}$ with $k$ terms. Given an operator $A$ on $\mathcal{H}$ we say that an operator $B$ on $\oplus^k \mathcal{H}$ is a *total* dilation of $A$, or that $B$ *totally* dilates $A$, if we can write

$$B = \begin{pmatrix} A & \star & \cdots \\ \star & A & \ddots \\ \vdots & \ddots & \ddots \end{pmatrix},$$

that is if the operator diagonal of $B$ consists of a repetition of $A$. Clearly this notion has also a natural extension when $A$ acts on any space $\mathcal{F}$ with $\dim \mathcal{F} = \dim \mathcal{H}$. Let $\{A_j\}_{j=0}^n$ be a family of operators on $\mathcal{H}$ and let $\{B_j\}_{j=0}^n$ be a family of operators on $\oplus^k \mathcal{H}$. We say that $\{B_j\}_{j=0}^n$ *totally* dilates $\{A_j\}_{j=0}^n$ if we can write, with respect to a (hilbertian) basis of $\mathcal{H}$,

$$B_0 = \begin{pmatrix} A_0 & \star & \cdots \\ \star & A_0 & \ddots \\ \vdots & \ddots & \ddots \end{pmatrix}, \quad \ldots, \quad B_n = \begin{pmatrix} A_n & \star & \cdots \\ \star & A_n & \ddots \\ \vdots & \ddots & \ddots \end{pmatrix}.$$

We give five examples of total dilations:

**Example 1.1.** A $2n \times 2n$ antisymmetric real matrix $A$ totally dilates the $n$-dimensional zero operator: with respect to a suitable decomposition

$$A = \begin{pmatrix} 0 & -B \\ B & 0 \end{pmatrix}$$

for some symmetric real $n$-by-$n$ matrix $B$.

**Example 1.2.** Any operator $A$ on $\mathcal{H}$ can be totally dilated into a normal operator $N$ on $\mathcal{H} \oplus \mathcal{H}$ by setting

$$N = \begin{pmatrix} A & A^* \\ A^* & A \end{pmatrix}.$$

**Example 1.3.** Denote by $\tau(A)$ the normalized trace $(1/n)\text{Tr}A$ of an operator $A$ on an $n$-dimensional space. Then the scalar $\tau(A)$ can be totally dilated into $A$. For an operator acting on a real space and for a hermitian operator the proof is easy. When $A$ is a general operator on a complex space, this result follows from the Hausdorff-Toeplitz Theorem (see [4, p. 20]).

**Example 1.4.** Any contraction $A$ on a finite dimensional space $\mathcal{H}$ can be totally dilated into a unitary operator $U$ on $\oplus^k \mathcal{H}$ for any integer $k \geq 2$. Indeed by considering the polar decomposition $A = V|A|$, it suffices to construct a total unitary dilation $W$ of $|A|$ and then to take $U = (\oplus^k V) \cdot W$. The construction of a total unitary dilation on $\oplus^k \mathcal{H}$ for a positive contraction $X$ on $\mathcal{H}$ is easy: Let $\{x_j\}_{j=1}^n$ be the eigenvalues of $X$ repeated according to their multiplicities and let $\{U_j\}_{j=1}^n$ be $k \times k$ unitary matrices such that $\tau(U_j) = x_j$. Example 1.3 and an obvious matrix manipulation show that $\oplus_{j=1}^n U_j$ totally dilates $X$.

**Example 1.5.** Let $\{A_k\}_{k=1}^n$ be a family of operators on $\mathcal{H}$ and let $\{B_k\}_{k=1}^n$ be the family of of operators acting on $\oplus^n \mathcal{H}$ defined by

$$B_k = \begin{pmatrix} A_k & A_{k-1} & \cdots \\ A_{k+1} & A_k & \cdots \\ \vdots & \vdots & \ddots \end{pmatrix}.$$

Then $\{B_k\}_{k=1}^n$ is a commuting family which totally dilates $\{A_k\}_{k=1}^n$. (We set $A_0 = A_n$, $A_{-1} = A_{n-1}, \ldots$)

In the last example above, the dilations do not preserve properties such as positivity, self-adjointness or normality. Using larger dilations we may preserve these properties:

**Proposition 1.6.** *Let $\{A_j\}_{j=0}^n$ be operators on a space $\mathcal{H}$. Then there exist operators $\{B_j\}_{j=0}^n$ on $\oplus^k \mathcal{H}$, where $k = 2^n$, such that*
(a)   *For $i \neq j$, $B_i B_j = 0$.*
(b)   *$\{B_j\}_{j=0}^n$ totally dilates $\{A_j\}_{j=0}^n$.*
(c)   *If the $A_j$'s are positive (resp. hermitian, normal) then the $B_j$'s are of the same type.*

**Proof.** Given a pair $A_0$, $A_1$ of operators, construct

$$S = \begin{pmatrix} A_0 & A_0 \\ A_0 & A_0 \end{pmatrix} \quad \text{and} \quad T = \begin{pmatrix} A_1 & -A_1 \\ -A_1 & A_1 \end{pmatrix}.$$

Then $ST = TS = 0$. We then proceed by induction. We have just proved the case of $n = 1$. Assume that the result holds for $n - 1$. Thus we have a family $\mathcal{C} = \{C_j\}_{j=0}^{n-1}$ which totally dilates $\{A_j\}_{j=0}^{n-1}$. Moreover $\mathcal{C}$ acts on a space $\mathcal{G}$, $\dim \mathcal{G} = 2^{n-1} \dim \mathcal{H}$. We dilate $A_n$ to an operator $C_n$ on $\mathcal{G}$ by setting $C_n = A_n \oplus \ldots \oplus A_n$, $2^{n-1}$ terms. We then consider the operators on $\mathcal{F} = \mathcal{G} \oplus \mathcal{G}$ defined by

$$B_j = \begin{pmatrix} C_j & C_j \\ C_j & C_j \end{pmatrix} \quad \text{for} \quad 0 \leq j < n \quad \text{and} \quad B_n = \begin{pmatrix} C_n & -C_n \\ -C_n & C_n \end{pmatrix}.$$

The family $\{B_j\}_{j=0}^n$ has the required properties.   ◇

If $\mathcal{H}$ is a space with an even finite dimension, we then say that the orthogonal decomposition $\mathcal{H} = \mathcal{F} \oplus \mathcal{F}^\perp$ is a *halving* decomposition whenever $\dim \mathcal{F} = (1/2) \dim \mathcal{H}$.

**Theorem 1.7.** *Let $A$ be an operator on a space $\mathcal{H}$ with an even finite dimension. Then there exists a halving decomposition $\mathcal{H} = \mathcal{F} \oplus \mathcal{F}^\perp$ for which we have a total dilation*

$$A = \begin{pmatrix} B & \star \\ \star & B \end{pmatrix}.$$

**Proof.** Choose a halving decomposition of $\mathcal{H}$ for which we have a matrix representation of $\text{Re}\,A$ of the following form

$$\text{Re}\,A = \begin{pmatrix} S & 0 \\ 0 & T \end{pmatrix}.$$

Consequently with respect to this decomposition we must have

$$A = \begin{pmatrix} Y & X \\ -X^* & Z \end{pmatrix}.$$

Let $X = U|X|$ and $Y_0 = U^*YU$. We have

$$\begin{pmatrix} U^* & 0 \\ 0 & I \end{pmatrix} A \begin{pmatrix} U & 0 \\ 0 & I \end{pmatrix} = \begin{pmatrix} U^* & 0 \\ 0 & I \end{pmatrix} \begin{pmatrix} Y & U|X| \\ -|X|U^* & Z \end{pmatrix} \begin{pmatrix} U & 0 \\ 0 & I \end{pmatrix}$$
$$= \begin{pmatrix} Y_0 & |X| \\ -|X| & Z \end{pmatrix}.$$

Now observe that

$$\frac{1}{\sqrt{2}} \begin{pmatrix} I & -I \\ I & I \end{pmatrix} \begin{pmatrix} Y_0 & |X| \\ -|X| & Z \end{pmatrix} \frac{1}{\sqrt{2}} \begin{pmatrix} I & I \\ -I & I \end{pmatrix}$$
$$= \begin{pmatrix} (Y_0 + Z)/2 & \star \\ \star & (Y_0 + Z)/2 \end{pmatrix}.$$

Thus, using two unitary congruence we have exhibited an operator totally dilated into $A$. $\diamondsuit$

**Remark 1.8.** The proof of Theorem 1.7 is easy for a normal operator $A$: consider a representation $A = \begin{pmatrix} S & 0 \\ 0 & T \end{pmatrix}$ and use the unitary conjugation by $\frac{1}{\sqrt{2}} \begin{pmatrix} I & I \\ -I & I \end{pmatrix}$. Applying this to $X^*X$, for an operator $X$ on an even dimensional space, we note that there exists a halving projection $E$ such that $XE$ and $XE^\perp$ have the same singular values (indeed $EX^*XE$ and $E^\perp X^*XE^\perp$ are unitarily equivalent).

**Problems 1.9.** (a) Does Theorem 1.7 extend to infinite dimensional spaces ? (b) Let $\mathcal{H}$, $\mathcal{F}$ be two finite dimensional spaces with $\dim \mathcal{H} = k \dim \mathcal{F}$ for an integer $k$. Is any operator $A$ on $\mathcal{H}$ a total dilation of some operator $B$ on $\mathcal{F}$ ?

The author has the feeling that the two questions above have a positive answer.

## 2. Constructions of monotone dilations

Recall that the notion of a monotone family of positive or hermitian operators has been discussed in the introduction.

*2.1 Algorithmic constructions of monotone dilations*

Given an operator $A$ on $\mathcal{H}$ and an integer $k > 0$ we define the following total dilations of $A$ on $\oplus^k \mathcal{H}$:

$$A(k) = \begin{pmatrix} A & 0 & \cdots \\ 0 & A & \ddots \\ \vdots & \ddots & \ddots \end{pmatrix} \quad \text{and} \quad A[k] = \begin{pmatrix} A & A & \cdots \\ A & A & \ddots \\ \vdots & \vdots & \ddots \end{pmatrix}.$$

Therefore, denoting by $I_k$ the $k$-by-$k$ identity matrix and by $E_k$ the $k$-by-$k$ matrix whose entries all equal to 1, we have $A(k) = A \otimes I_k$ and $A[k] = A \otimes E_k$. Note that $(1/k)E_k$ is a (rank one) projection, consequently, when $A$ is positive so is $A[k]$. For $k > 1$ we introduce another total dilation of $A$ on $\oplus^k \mathcal{H}$ by setting

$$A\langle k\rangle = \begin{pmatrix} A & \frac{I-A}{k-1} & \cdots \\ \frac{I-A}{k-1} & A & \ddots \\ \vdots & \ddots & \ddots \end{pmatrix}.$$

Thus we have

$$A\langle k\rangle = \left(\frac{I-A}{k-1}\right)[k] + \left(\frac{kA-I}{k-1}\right)(k).$$

If $A$ is a positive operator satisfying $I \geq A \geq (1/k)I$ the above relation shows that $A\langle k\rangle$ is a positive operator. Given two operators $A$, $B$ on $\mathcal{H}$ one can check that $A[k]$ and $B\langle k\rangle$ commute, in fact

$$A[k]B\langle k\rangle = A[k] = B\langle k\rangle A[k].$$

If both $A$ and $B$ are positive, a more precise result holds.

**Proposition 2.1.** *Let $(A, B)$ be a pair of positive operators on $\mathcal{H}$ and assume that $I \geq B \geq (1/k)I$ for some integer $k > 0$. Then $(A[k], B\langle k\rangle)$ is a monotone pair of positive operators which totally dilates $(A, B)$.*

Proposition 2.1 is just a restatement of Theorem 2.11 in [2]. The next result is a generalization for more general families than pairs. It is convenient to introduce some notations. First an expression like $A(k)\langle l\rangle[m]$ should be understood in the following way: begin by constructing $B = A(k)$, then construct $C = B\langle l\rangle$ and finally construct $C[m]$. Second, given a sequence $\{k_j\}_{j=1}^n$ of integers, we complete it with $k_{-1} = k_0 = k_{n+1} = 1$ and we set, for $0 \leq j \leq n$:

$$k'_j = \prod_{l=0}^{j-1} k_l \quad \text{and} \quad k''_j = \prod_{l=j+1}^{n} k_l \quad \text{(consequently } k'_0 = k''_n = 1\text{)}.$$

**Theorem 2.2.** *Let $\{A_j\}_{j=0}^n$ be positive operators on a space $\mathcal{H}$. Assume that for $j > 0$ we have integers $k_j > 0$ such that $I \geq A_j \geq (1/k_j)I$. Then there exist positive operators $\{B_j\}_{j=0}^n$ on $\oplus^k \mathcal{H}$, where $k = \prod_{j=1}^n k_j$, such that:*
(a) *$\{B_j\}_{j=0}^n$ is a monotone family of positive operators.*
(b) *$\{B_j\}_{j=0}^n$ totally dilates $\{A_j\}_{j=0}^n$.*
*A suitable choice for each $B_j$ is $A(k'_j)\langle k_j\rangle[k''_j]$.*

Multiplying by appropriate scalars, we note that the assumptions $I \geq A_j \geq (1/k_j)I$ may be replaced by $\mathrm{cond}(A_j) = ||A_j||.||A_j^{-1}|| \leq k_j \ (j > 0)$.

**Proof.** We proceed by induction. For $n = 1$, this is Theorem 2.11 in [2]. Assume that the result holds for $n - 1$. Let $\mathcal{A}_0 = \{A_j\}_{j=0}^{n-1}$. By the induction assumption there is a monotone family $\mathcal{C} = \{C_j\}_{j=0}^{n-1}$ which totally dilates $\mathcal{A}_0$. Furthermore $\mathcal{C}$ acts on a space $\mathcal{G}$ with $\dim \mathcal{G} = \prod_{j=1}^{n-1} k_j \dim \mathcal{H} = k'_n \dim \mathcal{H}$. Next, we dilate $A_n$ into an operator $C_n$ on $\mathcal{G}$ by setting $C_n = A_n(k'_n)$. To prove the theorem it now suffices to show that we can totally dilate the family $\mathcal{C}' = \{C_j\}_{j=0}^n$ on $\mathcal{G}$ into a monotone family $\mathcal{B} = \{B_j\}_{j=0}^n$ on a larger space $\mathcal{F}$ with $\dim \mathcal{F} = k_n \dim \mathcal{G}$.

To this purpose we consider on $\mathcal{F} = \mathcal{G} \oplus \ldots \oplus \mathcal{G}$, $k_n$ terms, the following operators: for $0 \leq j \leq n-1$,

$$B_j = \begin{pmatrix} C_j & C_j & \cdots \\ C_j & C_j & \cdots \\ \vdots & \vdots & \ddots \end{pmatrix}$$

and for $j = n$

$$B_n = \begin{pmatrix} C_n & \frac{I-C_n}{k_n - 1} & \cdots \\ \frac{I-C_n}{k_n - 1} & C_n & \ddots \\ \vdots & \ddots & \ddots \end{pmatrix}.$$

Because $\{C_j\}_{j=0}^{n-1}$ is a monotone family, so is $\{B_j\}_{j=0}^{n-1}$ (recall that $B_j = C_j \otimes E_{k_n}$ for $j < n$ where $E_p$ is, up to a scalar multiple, a rank one projection). Reasoning as in the proof of Theorem 2.11 [2] we obtain that $(B_j, B_n)$, $0 \leq j < n$, are monotone pairs. Consequently $\{B_j\}_{j=0}^n$ is a monotone family (if $\{X_j\}_{j=0}^{n-1}$ is a monotone family and $(X_j, X_n)$ are monotone pairs, $j < n$, then $\{X_j\}_{j=0}^n$ is a monotone family). Finally a close look to our constructions reveals that the $B_j$'s are given by the formulae of the last part of the theorem. $\diamond$

**Corollary 2.3.** *Let $\{A_j\}_{j=0}^n$ be hermitian operators on a space $\mathcal{H}$. Then we can totally dilate them into a monotone family of hermitian operators on a larger space $\mathcal{F}$ with $\dim \mathcal{F} = 2^n \dim \mathcal{H}$.*

**Proof.** We set $A'_j = \alpha_j A_j + (3/4)I$ where $\alpha_j > 0$ is sufficiently small to have $1/2 I \leq A'_j \leq I$. We apply Theorem 2.2 to dilate $A'_j$ to $B'_j$. The operators $B_j = (1/\alpha_j) B'_j - (3/4\alpha_j) I$ are the wanted dilations. $\diamond$

We may note that the proofs of the two preceding results have an algorithmic nature. More precisely, let us consider a sequence of hermitians $\{A_j\}_{j=0}^n$. The Frobenius norm $\|A_j\|_2$ is easily computed. Setting $\alpha_j = 1/4\|A_j\|_2$ and applying Theorem 2.2 as in the proof of Corollary 2.3 we may easily construct a monotone family totally dilating $\{A_j\}_{j=0}^n$.

**Remark 2.4.** If $A$, $B$ are positive noninvertible operators, it is not possible, in general, to dilate them into a positive monotone pair. Let

$$A = \begin{pmatrix} 1 & 0 \\ 0 & 0 \end{pmatrix} \quad \text{and} \quad B = \begin{pmatrix} 0 & 0 \\ 0 & 1 \end{pmatrix}.$$

Suppose that $(S, T)$ is a positive, monotone dilation of $(A, B)$. We should have the matrix representations repectively to a basis $(e_1, \ldots, e_n)$ of some space

$$S = \begin{pmatrix} 1 & 0 & \star & \ldots \\ 0 & 0 & 0 & \ldots \\ \star & 0 & \star & \ldots \\ \vdots & \vdots & \vdots & \ddots \end{pmatrix}, \quad T = \begin{pmatrix} 0 & 0 & 0 & \ldots \\ 0 & 1 & \star & \ldots \\ 0 & \star & \star & \ldots \\ \vdots & \vdots & \vdots & \ddots \end{pmatrix}.$$

Since $(S, T)$ is supposed to be positive, monotone we would have one of the following relations: $\ker S \subset \ker T$ or $\ker T \subset \ker S$. Say $\ker S \subset \ker T$, we would deduce that $Te_1 = Te_2 = 0$ and we would reach a contradiction.

*2.2. Theoretical constructions of monotone dilations*

In the previous subsection we have constructed monotone dilations in a rather explicit way by using matrix manipulations. Now we give more theoretical constructions; the resulting dilations will act on more economical spaces but will not be total dilations. Our first construction uses a standard dilation argument in connection with the numerical range of an operator and we refer the reader to chapter 1 of [4] for a detailed discussion of the numerical range.

**Proposition 2.5** *Let $A$, $B$ be two strictly positive operators on a space $\mathcal{H}$. Then we can dilate them into a monotone pair of strictly positive operators on a larger space $\mathcal{F}$ with $\dim \mathcal{F} = 6 \dim \mathcal{H}$.*

**Proof.** Invertibility of $A$ and $B$ ensures the existence of a real $r > 0$ such that

$$S = \begin{pmatrix} A & A - rI \\ A - rI & A \end{pmatrix} \quad \text{and} \quad T = \begin{pmatrix} B & -B + rI \\ -B + rI & B \end{pmatrix}$$

are strictly positive operators. Moreover $ST = TS$. Hence $N = S + iT$ is a normal operator acting on $\mathcal{G} = \mathcal{H} \oplus \mathcal{H}$. Because $S > 0$ and $T > 0$, the spectrum of $N$, $\mathrm{Sp}N$, lies in the open quadrant of $\mathbf{C}$,

$$Q = \{\, z = x + iy \mid x > 0 \text{ and } y > 0 \,\}.$$

We may then find a triangle $\Delta = \{x_1 + iy_1, x_2 + iy_2, x_3 + iy_3\}$ in $Q$ such that

$$x_1 < x_2 < x_3 \quad \text{and} \quad y_1 < y_2 < y_3 \quad (*)$$

and $\mathrm{conv}\Delta \supset \mathrm{Sp}N$. A standard dilation argument shows that there is a normal operator $M$ acting on a space $\mathcal{F} \supset \mathcal{G}$, $\dim \mathcal{F} = 3 \dim \mathcal{G}$, such that $\mathrm{Sp}M = \Delta$ and $M_\mathcal{G} = N$. Therefore

$$(\mathrm{Re}M)_\mathcal{H} = (\mathrm{Re}N)_\mathcal{H} = A \quad \text{and} \quad (\mathrm{Im}M)_\mathcal{H} = (\mathrm{Im}N)_\mathcal{H} = B.$$

¿From $(*)$ we deduce that $(\mathrm{Re}M, \mathrm{Im}M)$ is a monotone pair dilating $(A, B)$. $\diamondsuit$

At a time when it was not so clear to the author that a sequence of $n + 1$ hermitians could be dilated into a commuting family, T. Ando has pointed out to the author [1] the fact that it was a straightforward consequence of Naimark's Dilation Theorem. More precisely this theorem entails that the multiplicative constant $2^n$ in Proposition 1.6 can be replaced, in case of positive or hermitian operators, by $n + 2$ (but then the dilations are no longer total). We refer the reader to [3, p. 260] for a modern proof of Naimark's Theorem. Here the only thing we would need to know is the following particular case: Given positive operators $\{A_j\}_{j=0}^n$ on $\mathcal{H}$ satisfying $\sum A_j = I$, we can dilate them into a family $\{Q_j\}_{j=0}^n$ of mutually orthogonal projections on a larger space $\mathcal{F} = \mathcal{G} \otimes \mathcal{H}$ in which $\dim \mathcal{G} = n + 2$. Actually, rather than Naimark's Theorem, we only need the following much more elementary statement. Let us say that an operator $B$ *essentially* acts on a subspace $\mathcal{E}$ if both the range and the corange of $B$ are contained in $\mathcal{E}$ (equivalently, $\text{ran} B \subset \mathcal{E}$ and $(\ker B)^\perp \subset \mathcal{E}$).

**Lemma 2.6.** *Fix an integer $n$ and a space $\mathcal{H}$. Then there exist a larger space $\mathcal{F}$, $\dim \mathcal{F} = (n+1) \dim \mathcal{H}$, and an orthogonal decomposition $\mathcal{F} = \mathcal{E}_0 \oplus \ldots \oplus \mathcal{E}_n$, in which $\dim \mathcal{E}_j = \dim \mathcal{H}$ for each $j$, such that: for every family of operators $\{A_j\}_{j=0}^n$ on $\mathcal{H}$ there is a family $\{B_j\}_{j=0}^n$ of operators on $\mathcal{F}$ with $B_j$ essentially acting on $\mathcal{E}_j$ and $A_j = (B_j)_\mathcal{H}$, $0 \leq j \leq n$. Moreover when the $A_j$'s are hermitian or positive, the $B_j$'s can be taken of the same type.*

Let us sketch the elementary proof of this lemma. First, choose subspaces $\{\mathcal{E}_j\}_{j=0}^n$ of $\mathcal{F} = \oplus^{n+1}\mathcal{H}$ in such a way that for each $j$ (a) $\dim \mathcal{E}_j = \dim \mathcal{H}$, (b) The projection $E_j$ from $\mathcal{F}$ onto $\mathcal{E}_j$ verifies: $(E_j)_\mathcal{H}$ is a strictly positive operator on $\mathcal{H}$. Now, fix an integer $j$ and observe that any vector $h \in \mathcal{H}$ can be lifted to a unique vector $h_j \in \mathcal{E}_j$ such that $Hh_j = h$, where $H$ is the projection onto $\mathcal{H}$. Consequently any rank one operator of the form $R = h \otimes h$, $h \in \mathcal{H}$, can be lifted into a positive rank one operator $T$ essentially acting on $\mathcal{E}_j$ such that $T_\mathcal{H} = R$. This ensures that given a general (resp. hermitian, positive) operator $A$ on $\mathcal{H}$ there exists a general (resp. hermitian, positive) operator $B$ essentially acting on $\mathcal{E}_j$ such that $B_\mathcal{H} = A$.

**Theorem 2.7.** *Let $\{A_j\}_{j=0}^n$ be hermitian operators on a space $\mathcal{H}$. Then we can dilate them into a monotone family of hermitian operators on a larger space $\mathcal{F}$ with $\dim \mathcal{F} = 2(n+1) \dim \mathcal{H} - 1$.*

**Proof.** By Lemma 2.6 we may dilate $\{A_j\}_{j=0}^n$ into a commuting family of hermitians $\{S_j\}_{j=0}^n$ on a larger space $\mathcal{G}$ with $\dim \mathcal{G} = (n+1)\dim \mathcal{H} = d$. Thus, there is a basis $\{g_k\}_{k=0}^d$ in $\mathcal{G}$ and real numbers $\{s_{j,k}\}$ such that

$$S_j = \sum_{k=0}^d s_{j,k} g_k \otimes g_k \qquad (0 \leq j \leq n).$$

We take for $\mathcal{F}$ a space of the form

$$\mathcal{F} = \mathcal{E}_0 \oplus \mathcal{E}_1 \oplus \ldots \oplus \mathcal{E}_d$$

in which $\dim \mathcal{E}_0 = 1$ and $g_0 \in \mathcal{E}_0$; and for $k > 0$, $\dim \mathcal{E}_k = 2$ and $g_k \in \mathcal{E}_k$. Hence, we have $\dim \mathcal{F} = 2(n+1)\dim \mathcal{H} - 1$.

For $k > 0$, let $\{e_{1,k}; e_{2,k}\}$ be a basis of $\mathcal{E}_k$ and suppose that $g_k = (e_{1,k} + e_{2,k})/\sqrt{2}$ (*). We set, for $0 \leq j \leq n$,

$$B_j = s_{j,0} g_0 \otimes g_0 + \sum_{k=1}^{d} (r_{j,k} e_{1,k} \otimes e_{1,k} + t_{j,k} e_{2,k} \otimes e_{2,k})$$

where the reals $r_{j,k}$ and $t_{j,k}$ are chosen in such a way that:

(1)  $s_{j,k} = (r_{j,k} + t_{j,k})/2, \qquad j = 0, \ldots n.$

(2)  $r_{j,d} < \ldots < r_{j,1} < s_{j,0} < t_{j,1} < \ldots < t_{j,d}, \qquad j = 0, \ldots n.$

¿From (1) and (*) we deduce that $S_j = (B_j)_\mathcal{G}$ so that $A_j = (B_j)_\mathcal{H}$.
¿From (2) we infer that $\{B_j\}_{j=0}^n$ is a monotone family. $\diamond$

We close this paper with the final observation:

**Remark 2.8.** The results of Section 2 still hold for infinite dimensional spaces (and then we simply have $\mathcal{F} = \mathcal{H} \bigoplus \mathcal{H}$). Also, we may consider real operators on real spaces as well as complex operators on complex spaces.

**References**
[1] T. Ando, private communication.
[2] J.-C. Bourin, Singular values of compressions, restrictions and dilations, Linear Algebra Appl., to appear (2002) .
[3] K.R. Davidson, $C^*$-algebras by example, Fieds Institute Monographs **6**, American Mathematical Society, 1996.
[4] R.A. Horn, C.R. Johnson, Topics in matrix analysis, Cambrige University Press, 1991.